\documentclass[12pt]{article}
\date{}

\usepackage[dvips]{graphicx}
\usepackage{fancyhdr}
\usepackage{fancybox,eufrak}

\newtheorem{theo}{Theorem}

\def \IR {{\rm I \! R}}

\def \1 {{\large {\rm 1 \! \! 1}}}
\def \( {{\large {\rm ( \! ( }}}
\def \) {{\large {\rm ) \! ) }}}

\newcommand{\cqfd}{\mbox{}\nolinebreak\hfill\rule{2mm}{2mm}}

\newcommand{\usr}[2]{{ \renewcommand{\arraystretch}{0.5} \begin{array}[t]{c}
#1 \\ \scriptstyle #2 \end{array} }}

\begin{document}
\title{A case study of an Hamilton-Jacobi equation by the Adomian
decompositional method}

\author{T. K. Edarh-Bossou \& B. Birregah \\
D\'epartement de Math\'ematiques, Universit\'e de Lom\'e, BP. 1515
Lom\'e Togo\\ \{ tedarh , bbirregah \}@tg.refer.org }

\maketitle
\begin{abstract}
\noindent {\it We present a study of the Adomian's Decomposition
Method (ADM) applied to the Hamilton-Jacobi equations $\; u_t\,
+\, H\left( u_x\right) =0$.\\
We recall the well known characteristics methods in the case of
this type of equations to justify the existence or not of
solutions. This yields that the ADM gives efficient solutions in
time only in $]0,\, T^*[$, where $T^*$ is the critical time of our
equation.}\\
\end{abstract}
\\


\section{ A presentation of the ADM}
We consider in this work the Hamilton-Jacobi equation
\begin{eqnarray}P_u\left\{
\begin{array}{lll}
\displaystyle u_t+H\left( u_x\right) \; & =\;0 & \quad x\in\IR,\, t>0\\
\displaystyle u(0,x)  \;& =\;  u_0(x)&  \quad x\in\IR
\label{e1}\end{array}\right.\end{eqnarray} where $u=u(x,t)$ and
$u_0=u_0(x)$ are functions on a Hilbert spaces, and
$H$ a non linear operator.\\
Let take $\displaystyle L=\frac{\partial }{\partial t}$ as the
linear operator. The principle of the ADM relies on the assumption
that the solution $u$ of equation (\ref{e1}) can be set as a serie
$\displaystyle u\; =\; \sum_{k=0}^{+ \infty} u_{k} $. In
\cite{birr} the authors considered the non linear part as a
function of $u_x$. This yields to:
$$\displaystyle H\left( \frac{\partial u}{\partial x} \right) = \sum_{k=0}^{+ \infty} A_{k}\left(
\frac{\partial u_0 }{\partial x},\, ... ,\, \frac{\partial u_k
}{\partial x} \right). $$ and
$$A_0 = H\left(\frac{\partial u_0}{\partial x} \right) \quad \hbox{ with }
u_0=u(0,x) $$
\begin{eqnarray}
\displaystyle A_{n+1} &  \displaystyle \; =\; \sum_{k=0}^{n}
(k+1)u_{k+1}^{,} \frac{\partial}{\partial u_k^{,}}A_{n} \qquad
\end{eqnarray}
where $\displaystyle u_k^{,} = \frac{\partial u_k}{\partial x}
$.\\

\noindent Several authors give some formulas on the Adomian's
polynomials $\displaystyle A_{k} $ in the case of an operator
depending exclusively on $u$ (see for example \cite{abb,
babolian}). For more related works in this field one can see
\cite{adomian}, \cite{mustafa}, \cite{birregah}, \cite{kaya1}, \cite{adom} and \cite{Cherr}.\\
We present here other formulas more suited to our analysis.

\begin{theo}
Let's consider that the operator $H$ of (\ref{e1}) is indefinitely
differentiable in Frechet sense .\\
We have:\\
\noindent
 $\begin{array}{ll}
 \displaystyle
A_0 =&  H\left( u_0^,\right) \quad \hbox{and }  \displaystyle A_n
= \displaystyle \sum_{k=1}^n \frac{1}{k!}\;
\sum_{p_1+p_2+...+p_k=n} H^{(k)}\left( u_0^, \right)\left(  u_{p_1
}^,\, , u_{p_2}^,\, , ...,\, u_{p_k}^,\right) \label{e3}
\end{array}$\\
where $H^{(k)}\left( u_0^,\right)$ is the $k$-th derivative at
$u_o^,$, and $\displaystyle u_k^,=\frac{\partial u_k}{\partial
x}$. \label{theo1}
\end{theo}

\noindent
{\bf Proof:}\\
We note $\displaystyle \psi (\lambda )=\sum_{i=1}^n \lambda ^i
\frac{\partial u_i}{\partial x}$. Thus\\
$\displaystyle A_n\, =\, \frac{1}{n!}\, \frac{d^n}{d\lambda^n}
H(\psi (\lambda ))\, =\, \displaystyle \frac{1}{n!}\,
\frac{d^n}{d\lambda^n}[H\circ \psi ]
(\lambda )$,\\
but $\displaystyle [H\circ \psi ](\lambda )= [H\circ \psi
](\lambda _0)+ \sum_{k=1}^{n}\frac{1}{k!}[H\circ \psi
]^{(k)}\left( \lambda _0\right) \underbrace{ \left( \lambda
-\lambda _0,\, \lambda -\lambda _0,\, .....,\,
\lambda -\lambda _0 \right)}_{\hbox{k times}}$\\
then we deduce that:\\

\noindent$
\begin{array}{ll}
\displaystyle H( \psi (\lambda )) =H( \psi (\lambda _0)) & \\
\displaystyle + \sum_{k=1}^{n}\frac{1}{k!}\; \sum_{p_1+p_2+\,
...\, + p_k=n}\; \frac{1}{p_1!\, p_2!\, ...\, p_k!} H^{(k)}\left(
\psi \left( \lambda _0 \right) \right)\left( \psi ^{(p_1)}(\lambda
_0)(\lambda -\lambda _0)^{p_1} ,\, .... ,\, \psi ^{(p_k)}(\lambda
_0)(\lambda -\lambda _0)^{p_k} \right) &
\end{array}
$\\

\noindent
For $\lambda _0 =0$, we have:\\
\medskip
$
\begin{array}{ll}
\displaystyle H( \psi (\lambda )) =H\left( \frac{\partial u_0}
{\partial x}\right) & \\
\displaystyle + \sum_{k=1}^{n}\frac{1}{k!}\; \sum_{p_1+p_2+\,
...\, + p_k=n}\; \frac{1}{p_1!\, p_2!\, ...\, p_k!} H^{(k)}\left(
\frac{\partial u_0} {\partial x}\right) \left( \lambda ^{p_1}
p_1!\frac{\partial u_{p_1}}{\partial x},\; ....,\; \lambda
^{p_k}p_k!\frac{\partial u_{p_k}}{\partial x}\right)
 & \\
\displaystyle =H\left( \frac{\partial u_0}
{\partial x}\right) & \\
\displaystyle + \sum_{k=1}^{n}\frac{1}{k!}\; \sum_{p_1+p_2+\,
...\, + p_k=n}\; \frac{p_1!\, p_2!\, ...\, p_k!} {p_1!\, p_2!\,
...\, p_k!} H^{(k)}\left( \frac{\partial u_0} {\partial x}\right)
\left( \frac{\partial u_{p_1}}{\partial x},\; ....,\;
\frac{\partial u_{p_k}}{\partial x}\right)\; \lambda ^{\overbrace{
p_1+\, ....\, +p_k}^{=n}} & \\
\displaystyle =H\left( \frac{\partial u_0}
{\partial x}\right) & \\
\displaystyle + \lambda ^n \sum_{k=1}^{n}\frac{1}{k!}\;
\sum_{p_1+p_2+\, ...\, + p_k=n}\; H^{(k)}\left( \frac{\partial
u_0} {\partial x}\right) \left( \frac{\partial u_{p_1}}{\partial
x},\; ....,\; \frac{\partial u_{p_k}}{\partial x}\right) &
\end{array}
$\\[0.5cm]

\noindent Thus: $\displaystyle \frac{d^n}{d\lambda ^n} H(\psi
(\lambda ))=n! \sum_{k=1}^{n}\frac{1}{k!}\; \sum_{p_1+p_2+\, ...\,
+ p_k=n}\; H^{(k)}\left( \frac{\partial u_0} {\partial x}\right)
\left( \frac{\partial u_{p_1}}{\partial x},\; ....,\;
\frac{\partial u_{p_k}}{\partial x}\right)
$\\[0.5cm]
This ends our demonstration. \cqfd

\begin{theo}
For the problem $P_u$, we have $\displaystyle A_{n-1}=\tilde
u_n(x)\, \frac{t^{n-1}} {(n-1)!}\; $ \\
and $\; \displaystyle u_n(x,t)=\tilde u_n(x)\, \frac{t^n}{n!}
\quad \forall n\ge 1$, with $$\displaystyle  \tilde u_{n}\,=\,
\sum_{k=1}^{n-1} \frac{1}{k!}\; \sum_{p_1+p_2+...+p_k=n-1}
\frac{n!}{p_1!p_2!\, ...\, p_k!}\, H^{(k)}\left( u_0^,
\right)\left( \tilde u_{p_1 }^,\, , \tilde u_{p_2}^,\, , ...,\,
\tilde u_{p_k}^,\right) $$ \label{theo2}
\end{theo}

\noindent
{\bf Proof:}\\
This will be done by recurrence.\\
In fact, $A_0=H(u_0^,)$,\\
then $\displaystyle u_1=\int^t_0 H(u_0^,) ds=H(u_0^,)t$, and
$\tilde u_1(x)=
H(u_0(x))$\\
Let's assume that $\displaystyle \forall k\le n \; \,
A_{k-1}=\tilde u_k(x)\, \frac{t^{k-1}}
{(k-1)!}\; $ and $\; \displaystyle u_k(x,t)=\tilde u_k(x)\, \frac{t^k}{k!}.$\\
From Theorem 1, we have: $\displaystyle A_n=\sum_{k=1}^n
\frac{1}{k!}\; \sum_{p_1+p_2+...+p_k=n} H^{(k)}\left( u_0^,
\right)\left(  u_{p_1 }^,\, , u_{p_2}^,\, , ...,\, u_{p_k}^,\right)$.\\
Then:\\
$
\begin{array}{ll}
\displaystyle A_n & \displaystyle = \sum_{k=1}^n \frac{1}{k!}\;
\sum_{p_1+p_2+...+p_k=n} H^{(k)}\left( u_0^,
\right)\left(  u_{p_1 }^,\, , u_{p_2}^,\, , ...,\, u_{p_k}^,\right)\\
 &\displaystyle =
 \sum_{k=1}^n \frac{1}{k!}\; \sum_{p_1+p_2+...+p_k=n} H^{(k)}\left( u_0^,
\right)\left( \tilde u_{p_1 }^,\, , \tilde u_{p_2}^,\, , ...,\,
\tilde u_{p_k}^,\right)\times \frac{t^{p_1}}{p_1!}\frac{t^{p_2}}{p_2!}\, ...
\, \frac{t^{p_k}}{p_k!}\\
 &\displaystyle = \frac{t^n}{n!}\,
 \sum_{k=1}^n \frac{1}{k!}\; \sum_{p_1+p_2+...+p_k=n}
\frac{n!}{p_1!p_2!\, ...\, p_k!}\, H^{(k)}\left( u_0^,
\right)\left( \tilde u_{p_1 }^,\, , \tilde u_{p_2}^,\, , ...,\,
\tilde u_{p_k}^,\right)\\
\end{array}
$\\

\noindent We have:
$$\displaystyle  \tilde u_{n+1}\,=\,
\sum_{k=1}^n \frac{1}{k!}\; \sum_{p_1+p_2+...+p_k=n}
\frac{n!}{p_1!p_2!\, ...\, p_k!}\, H^{(k)}\left( u_0^,
\right)\left( \tilde u_{p_1 }^,\, , \tilde u_{p_2}^,\, , ...,\,
\tilde u_{p_k}^,\right) \; ,$$

\noindent and, $ \displaystyle u_{n+1}\, =\int^t_0 A_n ds\, =
\tilde u_{n+1}(x)
\int^t_0\frac{s^n}{n!}ds$\\
this leads to the result. \cqfd\\

\noindent These two formulas show that $u_n$ closely depends on
the derivatives $\displaystyle H^{(k)}\left(\frac{\partial u_0
}{\partial x} \right)$, consequently on
$\displaystyle\frac{\partial u_0}{\partial x} $. On the other side
the form of $\displaystyle u_n(x,t)=\tilde u_n(x)\,
\frac{t^n}{n!}$ can help us to set some conjectures about the
convergence of the serie $\displaystyle \sum_{k=0}^{\infty } u_k(x,t)$.\\
We have:
$\displaystyle \left| \frac{u_{n+1}}{u_{n}} \right| =\, \left| \frac{\tilde u_{n+1(x)}}{\tilde u_{n}(x)} \right|
\, \frac{\left| t\right|}{n+1} \le \, \left| \frac{\tilde u_{n+1(x)}}{\tilde u_{n}(x)} \right| \, |t|$\\

\noindent Using the well known rule of d'Alembert,  one can see
that the convergence of $\displaystyle \sum_{k=0}^{\infty }
u_k(x,t)$ depends on $\displaystyle H^{(k)}\left(\frac{\partial
}{\partial x} u_0
\right)$ and $\displaystyle\frac{\partial }{\partial x} u_0$ .\\

\section{ Characteristics method}

\noindent We still consider the problem $P_u$ in (\ref{e1}).\\
By deriving formally in space ($P_u$) and taking $\displaystyle v=
\frac{\partial u}{\partial x}$, we have the conservation law
equation
\begin{eqnarray}
v_t\, +\, \frac{\partial }{\partial x}\left[ H\left(v
\right)\right] =0 \label{ee}
\end{eqnarray}
with the initial condition $\displaystyle v(x,0)=v_0(x)= \frac{\partial u_0}{\partial x}$.\\
Assuming that $v$ is derivable in space, we then have:
\begin{eqnarray}
v_t\, +\, H'\left(v \right)\, v_x =0 \label{e4} \end{eqnarray}
\noindent Using $\displaystyle a(v)=H'(v)$, we have: $v_t\, +\,
a\left(v \right)\, v_x =0$\\ Let's consider in the plan $(x,t)$
the characteristic strips $\displaystyle x \mapsto x(t)$ solutions
of the equation
\begin{eqnarray} \left\{
\begin{array}{l}
\displaystyle
\frac{dx}{dt} = a\left( v\left( x(t),\, t\right) \right) \label{e5}\\
x(0)=x_0
\end{array}\right. \end{eqnarray}
Writing $\displaystyle w(t)= v\left( x(t),\, t\right)$, we have:
\begin{eqnarray}
\frac{d}{dt} w(t)\,=\, \frac{d}{dt} v\left( (x,t),\, t\right) &
\displaystyle = \left( v_t +  \frac{dx}{dt} \frac{ \partial v}{
\partial x}\right) \left( (x,t),\, t \right) \nonumber \\
& = \displaystyle \left( v_t +  a \frac{ \partial v}{
\partial x}\right) \left( (x,t),\, t \right) = \, 0\nonumber
\end{eqnarray}

\noindent This leads to: $\displaystyle w(t)= w_0$; i.e. $v$ is constant  on the the characteristic strips.\\
\noindent From (\ref{e5}) we have
$\displaystyle \frac{dx}{dt}\, = \, a\left( v_0\left(  x_0 \right)
\right) \nonumber $,
then the characteristic strip from $x_0$ is a straight line with
equation given by:
\begin{eqnarray}
\displaystyle x(t)= a\left( v_0\left(  x_0 \right) \right)\, t +
x_0
\end{eqnarray}

{\setlength{\unitlength}{8mm}
\begin{figure}[h]
\begin{center}
\begin{minipage}[c]{0.40\linewidth}
\begin{center}
\begin{picture}(8,8)
\includegraphics[angle=0, width=6cm, height=6cm]{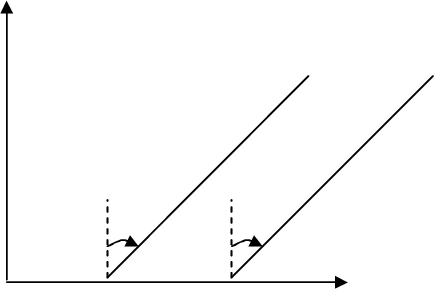}
\put(-3, 1.5){{\footnotesize $a\left(
v_0\left(x_2\right)\right)$}} \put(-6, 1.5){{\footnotesize
$a\left( v_0\left(x_1\right)\right)$}} \put(-6,-0.5){$x_1$}
\put(-4,-0.5){$x_2$} \put(-1,0){$x$} \put(-7.8,7){$t$}
\end{picture}
\end{center}
\caption{\footnotesize {\it Characteristics with non decreasing
$a\left( v_0\left( . \right) \right)$}}
 \label{caract11}
\end{minipage}%
\hfill
\begin{minipage}[c]{0.40\linewidth}
\begin{center}
\begin{picture}(8,8)
\includegraphics[angle=0, width=6cm, height=6cm]{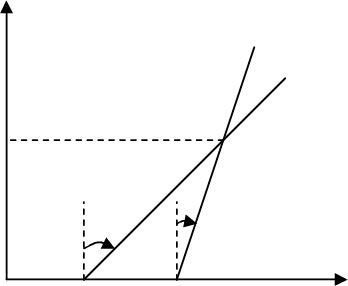}
\put(-3, 1.5){{\footnotesize $a\left(
v_0\left(x_2\right)\right)$}} \put(-6, 1.5){{\footnotesize
$a\left( v_0\left(x_1\right)\right)$}} \put(-6,-0.5){$x_1$}
\put(-4,-0.5){$x_2$} \put(0,0){$x$} \put(-7.8,7){$t$}
\put(-7.8,3.7){$\tilde T$} \put(-2.8,2.5){$ \leftarrow
u=u_0\left(x_1\right)$} \put(-1.8,4.3){$ \leftarrow u=u_0\left(
x_1\right)$}
\end{picture}
\end{center}
\caption{\footnotesize {\it Characteristics with decreasing
$a\left( v_0\left( . \right) \right)$ }}\label{caract1}
\end{minipage}
\end{center}
\label{neumann}
\end{figure}}
\noindent If $a\left( v_0\left(  . \right) \right)$ is non
decreasing, we obtain regular solutions , else one can observe
losts of regularity (see fig.\ref{caract11}), i.e. the solution
lost in regularity in time at $\tilde T$.

\noindent The time $\displaystyle T^* = \inf \{ \tilde T \}$ is
called critical time of the equation (\ref{ee}) and is calculated
by \cite{leroux, stanley, sethian87}:
\begin{eqnarray} T^* = - \frac{1}{\usr{\min }{x\in \IR }\left( \displaystyle H'\circ
\frac{\partial u_0}{\partial x}\right)'(x)} \,  =\, -
\frac{1}{\usr{\min }{x\in \IR }\left( a\left( v_0 \right)
\right)'(x) }. \label{critique}
\end{eqnarray}
 We can see that $\displaystyle H'\circ \frac{\partial u_0}{\partial
 x}$ is one of the terms $\displaystyle H^{(k)}\circ \frac{\partial u_0}{\partial
 x}$ in the Adomian polynomials $A_n$.\\
This leads us to say that the radius of convergence in time of the
series $\displaystyle \sum_{k=0}^{\infty } u_k(x,t)$ links closely
to the critical time $T^*$.

\section{Applications}

\noindent
\subsection{Example 1}
Let's consider the equation:
$$\left\{
\begin{array}{l}
\displaystyle
u_t + \frac{1}{2}\left( u_x \right)^2=0 \\[0.5cm]
u(x,0)=u_0(x)= - x^2
\end{array}
\right.$$ The associate conservation law problem is the well-known
Burger's equation \cite{JMburg}:
$$\left\{
\begin{array}{l}
\displaystyle
v_t +\left( \frac{v^2}{2} \right)_x=0 \\[0.5cm]
v(x,0)=v_0(x)= -2x
\end{array}
\right.$$ \noindent We then have $\displaystyle
H(v)=\frac{1}{2}v^2$. Thus $H'(v)=v\equiv a(v)$, and:
$$\displaystyle T^* \, =\, - \frac{1}{\usr{\min }{x\in \IR }\left(
a\left( v_0 \right) \right)'(x) } \, =\, - \frac{1}{\usr{\min
}{x\in \IR }\left( v_0 \right)'(x) }\, =\, \frac{1}{2} $$\\

\noindent Consider Theorem 1, the ADM gives:
$$\begin{array}{lcl}
u_0 &=& - x^2 \\
A_0 &=& \displaystyle \frac{1}{2} (u_o^{,})^2 \, =\, 2x^2 \\
u_1 &= & \displaystyle -\int A_0 dt = - 2tx^2 \\
A_1 &=& u_o^{,}u_1^{,} \\
u_2 &=& -4t^2x^2 \\
... & ... & \\
u_n &=& -x^2(2t)^n \\
\end{array}
$$
\noindent So: $$\displaystyle u(t,x) =-x^2 \sum_{n=0}^{\infty
}(2t)^n .$$ For $|t|< 1/2 $ we have:
$$\displaystyle u(t,x) = \frac{x^2}{1-2t}$$
\noindent  One can easily verify that after $t=1/2$,
$\displaystyle u(t,x) = \frac{x^2}{1-2t}$ remains a solution of
our equation. But this is not the physically true. As we can see,
$u(t,x)$ isn't regular at $t=1/2$. Recalling the fact that this
equation model frontier evolution (for example in mathematical
morphology \cite{article1, birregah}), we know that this
irregularity spreads forward for $t>1/2$. At this stage, only
generalized or entropic solutions survived after $\displaystyle
T^*=\frac{1}{2}$. We get these solutions by using Kruskov's
formula (see for example \cite{leroux}). Numerically these
solutions can be computed with Godunov or Hamilton-Jacobi
schemes.(see \cite{sethian87, sethian})

\subsection{Example 2}
$$\left\{
\begin{array}{l}
\displaystyle
u_t + \sqrt{1+ u_x ^2}=0 \\[0.5cm]
u(x,0)=u_0(x)
\end{array}
\right.$$
{\bf Case 1:} $u_0(x)=ax+b,\; a,b\in \IR $.\\
We can notice that the characteristic strips, in this case, are
parallel lines. So this equation has a regular solution at all
time. The viscosity solutions theory (see \cite{stromberg, ugv})
stipulates that the solutions are parallel lines to the initial
condition. One can verify then $H'\left( u_0'\right)$
is constant. Then $T^*=+\infty $.\\
With the ADM, we have:
$$\begin{array}{lcl}
u_o & = & ax+b \\
& & \\
A_0 & = & \sqrt{1+a^2} \\
 & & \\
u_1 & = & t\sqrt{1+a^2} \\
& & \\
A_1 & = & 0 \\
 & & \\
u_2 & = & 0 \\
... &... &... \\
 & & \\
u_n & = & 0  \qquad \forall n\ge 2
\end{array}$$

\noindent Thus:
$$u(t,x)= ax+b + t\sqrt{1+a^2} $$
This solution is exactly the viscosity solution given by the
well-known Lax formula \cite{stromberg, ugv}.\\

\noindent {\bf Case 2:} $u_0(x)=\sin x$

\noindent We then have $H(v)=\sqrt{1+v^2}$.\\
Thus $\displaystyle H'(v)=\frac{v}{\sqrt{1+v^2}} \equiv a(v)$,\\
and: \noindent $\displaystyle T^* \, =\, - \frac{1}{\usr{\min
}{x\in \IR }\left( a\left( v_0 \right) \right)'(x) } \,  =\, -
\frac{1}{\usr{\min
}{x\in \IR }\displaystyle \frac{d}{dx}\left( \frac{v_0}{\sqrt{1+v_0^2}} \right)(x) }$\\

\noindent $\displaystyle \frac{d}{dx}\left(
\frac{v_0}{\sqrt{1+v_0^2}} \right)(x) \, =\, \frac{d}{dx}\left(
\frac{\cos x}{\sqrt{1+\cos ^2 x }} \right)\, =\, \frac{-\sin
x}{\left( 1+ \cos ^2 x\right)
\sqrt{1+\cos ^2 x}}$\\

\noindent It is easy to see that:
$$\usr{\min}{x\in \IR } \, \frac{-\sin x}{\left( 1+ \cos ^2 x\right)
\sqrt{1+\cos ^2 x}}\, =\, \left[ \frac{-\sin x}{\left( 1+ \cos ^2
x\right) \sqrt{1+\cos ^2 x}}\right]_{x=\frac{\pi }{2}} .$$

\noindent Thus $\displaystyle T^* =\,1$.\\

\noindent
The ADM gives:\\
\noindent $\begin{array}{lcl}
& & \\
u_0 & = & \sin x \\
& & \\
u_1 & = & t\sqrt{1+\cos^2 x} \\
& & \\
u_2 & = & \displaystyle
-\frac{t^2}{2} \,\frac{\cos^2 x \sin x}{1+\cos^2 x} \\
& & \\
\end{array}$
\\
\noindent $\begin{array}{lcl} u_3 & = & \displaystyle
\frac{t^3}{6(1+\cos^2 x)^{7/2}} \, (2\cos^2 x \sin^2 x +\cos^4 x
\sin^2 x
-\cos^4 x -2\cos^6 x -\cos^8 x )\\
& & \\
u_4 & = & \displaystyle \frac{t^4}{24(1+\cos^2 x)^{5}} \,(10\cos^4
x \sin^2 x+13\cos^4 x \sin x-
4\cos^3 x \sin^2 x \\
 & &  \qquad \qquad \qquad + 24\cos^6 x \sin x+ 8\cos^8 x \sin x+
7\cos^{10} x\sin x-6\cos^2 x \sin^3 x \\
 & & \qquad \qquad \qquad -6\cos^4 x \sin^3 x+ 3\cos^2 x \sin x)\\
...&... & ... \\
\end{array}$\\
{\setlength{\unitlength}{8mm}
\begin{figure}[ht]
\begin{center}
\begin{minipage}[c]{0.40\linewidth}
\begin{center}
\begin{picture}(6.5,6.5)
\put(15,0.2){(b)} \put(4,0.2){(a)}
\includegraphics[angle=0, width=6cm, height=6cm]{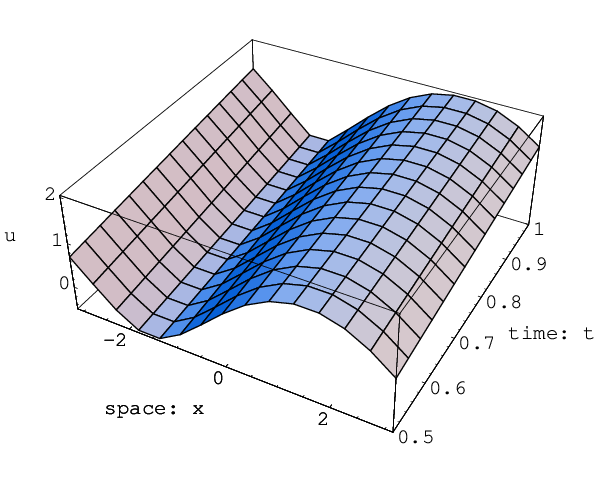}
\end{picture}
\end{center}
\end{minipage}%
\hfill
\begin{minipage}[c]{0.40\linewidth}
\begin{center}
\begin{picture}(6.5,6.5)
\includegraphics[angle=0, width=6cm, height=6cm]{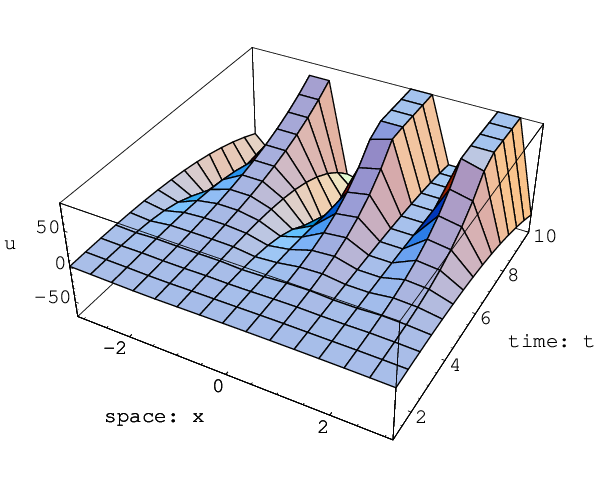}
\end{picture}
\end{center}
\end{minipage}%
\end{center}
\caption{ \footnotesize {\it $\displaystyle
\sum_{k=0}^{4}u_k(x,t)$ for (a) $0.5\le t\le 1 $ and (b) $1\le
t\le 10$ }} \label{3D}
\end{figure}}


{\setlength{\unitlength}{8mm}
\begin{figure}[ht]
\begin{center}
\begin{minipage}[c]{0.40\linewidth}
\begin{center}
\begin{picture}(6.5,6.5)
\put(16,0.5){(b)} \put(5,0.5){(a)}
\includegraphics[angle=0, width=6cm, height=6cm]{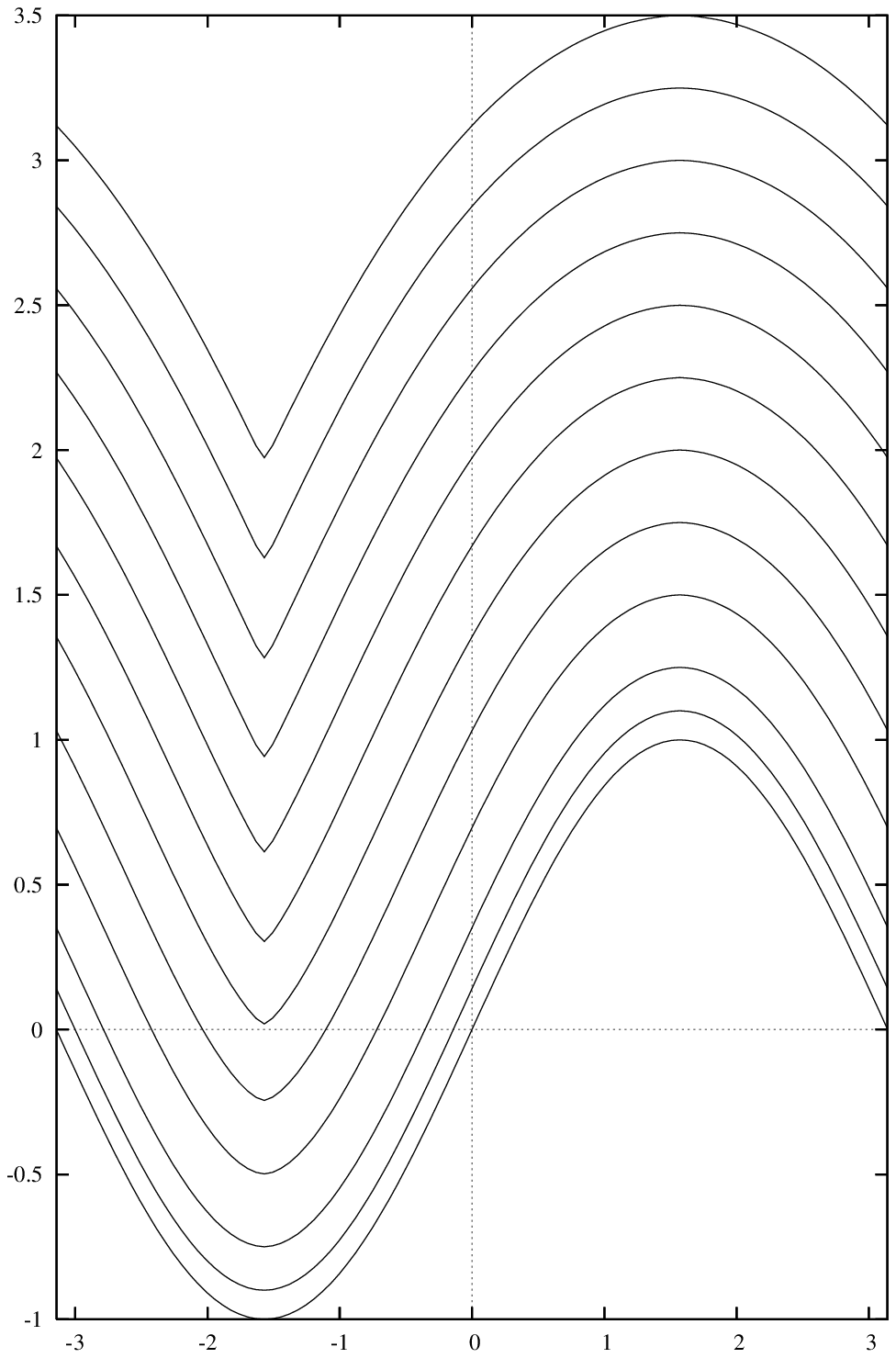}
\end{picture}
\end{center}
\end{minipage}%
\hfill
\begin{minipage}[c]{0.40\linewidth}
\begin{center}
\begin{picture}(6.5,6.5)
\includegraphics[angle=0, width=6cm, height=6cm]{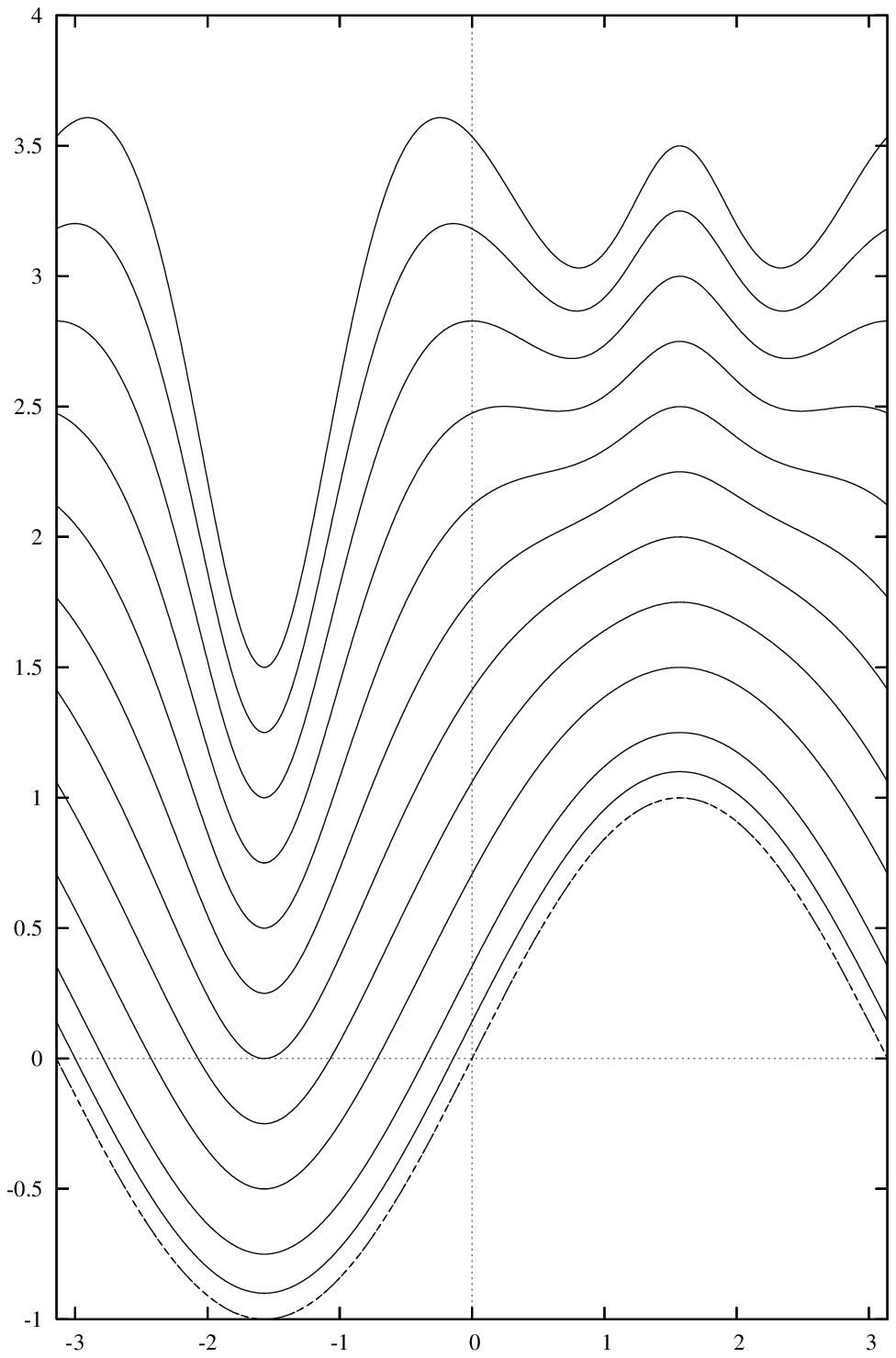}
\end{picture}
\end{center}
\end{minipage}%
\end{center}
\caption{ \footnotesize {\it (a) Numerical solution with
Hamilton-Jacobi scheme and (b) $\displaystyle
\sum_{k=0}^{4}u_k(x,t)$ at different times ( upwards : $t=0,\,
0.1,\, 0.25,\, 0.5,\, 0.75,\, 1,\, 1.25,\, 1.5,\, 1.75,\, 2,\,
2.25,\, 2.5 $). } }\label{1D}
\end{figure}}
\noindent One can see as shown on figures \ref{3D} and \ref{1D}
(a)-(b) that for $0\le t\le T^*$ (with $T^*=1$ the ADM gives an
accurate
approximation  of the solution.\\
\noindent This confirms that the ADM, for the case of
Hamilton-Jacobi equations, gives a good approximation of the
solution for only $t\le T^*$ even with few terms $u_n$.
Nevertheless at $t>\, T^*$
the computation of several terms of the series $\left( u_n\right)$ doesn't ameliorate the solution.\\

\section{Conclusion}

\noindent We show in this paper that, for Hamilton-Jacobi
equations and generally hyperbolic equations the ADM gives
accurate solution or a good approximation of the solution only for
$t\le \, T^*$. We deduce from there that ADM doesn't gives global
solution in time for this type of equation.\\

\end{document}